\documentclass[12pt]{article}

\usepackage{amsmath,amsfonts, amssymb}

\begin{document}
\bibliographystyle{plain}
\pagenumbering{arabic}
\raggedbottom

\newtheorem{theorem}{Theorem}%[section]
\newtheorem{lemma}[theorem]{Lemma}
\newtheorem{proposition}[theorem]{Proposition}
\newtheorem{corollary}[theorem]{Corollary}
\newtheorem{conjecture}[theorem]{Conjecture}
\newtheorem{definition}[theorem]{Definition}
\newtheorem{example}[theorem]{Example}
\newtheorem{condition}{Condition}
\newtheorem{main}{Theorem}
\renewcommand{\thefigure}{\thesection.\arabic{figure}}
\setlength{\parskip}{\parsep}
\setlength{\parindent}{0pt}
\setcounter{tocdepth}{1}

\def \outlineby #1#2#3{\vbox{\hrule\hbox{\vrule\kern #1% 
\vbox{\kern #2 #3\kern #2}\kern #1\vrule}\hrule}}%
\def \endbox {\outlineby{4pt}{4pt}{}}%

\newenvironment{proof}
{\noindent{\bf Proof\ }}{{\hfill \endbox
}\par\vskip2\parsep}

\hfuzz12pt

\newcommand{\cl}[1]{{\mathcal{C}}_{#1}}
\newcommand{\corr}{\mbox{Corr}}
\newcommand{\cov}{\mbox{Cov}}
\newcommand{\expo}{\mbox{Exp}}
\newcommand{\tr}[1]{\mbox{\rm{tr}}\left(#1 \right)}
\newcommand{\var}{\mbox{Var}}
\newcommand{\supp}{\mbox{supp}}
\newcommand{\tends}{\rightarrow \infty}
\newcommand{\ep}{{\mathbb {E}}}
\newcommand{\pr}{{\mathbb {P}}}
\newcommand{\re}{{\mathcal{R}}}
\newcommand{\rt}{\widetilde{\rho}}
\newcommand{\bds}{\begin{displaystyle}}
\newcommand{\eds}{\end{displaystyle}}
\newcommand{\vc}[1]{{\mathbf{#1}}}
\newcommand{\ra}[2]{{#1}^{(#2)}}
\newcommand{\vd}[1]{{\boldsymbol{#1}}}
\newcommand{\wir}{\widetilde{\rho}}

\newlength{\dsl}
\newcommand{\doublesub}[2]{\settowidth{\dsl}{$\scriptstyle
#2$}\parbox{\dsl}{\scriptsize\centering { \normalsize $\scriptstyle #1$}
\\ {\normalsize $\scriptstyle #2$}}}

\title{A conditional {Entropy Power Inequality} for dependent variables}
\author{Oliver Johnson}
\date{\today}
\maketitle

\makeatletter
\begin{abstract} 
We provide a condition under which a version of Shannon's Entropy Power
Inequality will hold for dependent variables. We provide information
inequalities extending those found in the independent case.
\renewcommand{\thefootnote}{}
\footnote{{\bf Addresses:} O.Johnson, Statistical Laboratory, 
Centre for Mathematical Sciences, Wilberforce Road, Cambridge, CB3 0WB,
UK.
Contact Email: {\tt otj1000@cam.ac.uk}. }
\renewcommand{\thefootnote}{\arabic{footnote}}
\setcounter{footnote}{0}
\makeatother
\end{abstract}

Shannon's Entropy Power Inequality states that 
\begin{theorem} \label{thm:epi} For independent random variables
$X,Y$ with densities, the entropy of the sum satisfies:
$$ 2^{2H(X+Y)} \geq 2^{2H(X)} + 2^{2H(Y)},$$
with equality if and only if $X,Y$ are normal. \end{theorem}
Apart from its intrinsic interest, it provides a sub-additive inequality for
sums of random variables and is thus an important part of the entropy-theoretic
proof of the Central Limit Theorem \cite{barron}. Whilst Shannon's proof
\cite{shannon} seems incomplete, in that he only checks that the necessary
conditions for a local maximum are satisfied, a rigorous proof is provided by
Stam \cite{stam} (see Blachman \cite{blachman}). This proof is based on a
related inequality concerning Fisher information:

\begin{lemma} \label{lem:fishinv}
For $X,Y$ with differentiable densities:
$$ \frac{1}{J(X+Y)} \geq \frac{1}{J(X)} + \frac{1}{J(Y)}.$$
\end{lemma}

The relationship between Theorem \ref{thm:epi} and 
Lemma \ref{lem:fishinv} comes via de Bruijn's identity, which expresses
entropy as an integral of Fisher informations.

Now, Takano \cite{takano} provided conditions on the random 
variables $X,Y$, such that Theorem \ref{thm:epi} would still hold for
weakly dependent variables. In contrast, we change the equation, replacing
entropies by conditional entropies, providing alternative conditions for this
related result. Our approach is again to develop a Fisher information 
inequality, and to use an integral form of that to deduce the full result.

We consider random variables $X,Y$ with joint density $p(x,y)$ and
marginal densities $p_X(x),p_Y(y)$. We need to refer to score functions and 
Fisher informations. Write $\rho_X(x) = p_X'(x)/p_X(x)$  and $\rho_Y(y) = 
p_Y'(x)/p_Y(x)$. We write $\ra{p}{1}(x,y)$ for $\partial p(x,y)/\partial x$,
and similarly for $\ra{p}{2}(x,y)$, and $\ra{\rho}{1}(x,y) = 
\ra{p}{1}(x,y)/p(x,y)$, and similarly for $\ra{\rho}{2}(x,y)$. Now, we can 
define $J(X) = \ep \rho_X(X)^2$ and
$J(Y) = \ep \rho_Y(Y)^2$ for the Fisher informations of $X$ and $Y$, and
$J_{XX} = \ep \ra{rho}{1}(X,Y)^2$, $J_{YY} = \ep \ra{\rho}{2}(X,Y)^2$,
$J_{XY} = \ep \ra{\rho}{1}(X,Y) \ra{\rho}{2}(X,Y)$, similarly. We will need 
to consider terms of the form:
$M_{a,b}(x,y) = a (\ra{\rho}{1}(x,y) - \rho_X(x)) + b (\ra{\rho}{2}(x,y) - 
\rho_Y(y))$.
\begin{lemma}[Takano]
As in the independent case, we can express the score function $\rho_W$
of the sum $W = X+Y$ as a conditional expectation of $\rho_X(X,Y)$. 
\label{lem:fishcondsum}
$$ \rho_W(w) = \ep ( \ra{\rho}{2}(X,Y) | X+Y = w) = \ep ( \ra{\rho}{2}(X,Y) 
| X+Y=w).$$
\end{lemma}
\begin{proof} 
Since
$W = X+Y$ has density $p_W(w) = \int p(x,w-x) dx = \int p(w-y,y) dy$, 
we know that:
$$ \rho_Z(z) = \frac{p_W'(z)}{p_W(z)} = \int \frac{ \ra{p}{1}(w-y,y)}{p_W(w)} 
dy
= \int \frac{ \ra{p}{1}(w-y,y)}{p(w-y,y)} \frac{p(w-y,y)}{p_W(w)} dy,$$
hence the result follows. \end{proof} 
Using this, we establish the following proposition, the equivalent of Lemma 
\ref{lem:fishinv} for dependent variables, and which reduces to 
Lemma \ref{lem:fishinv} in the independent case:
\begin{proposition} \label{prop:fishdec} For random variables $X$,$Y$ with
differentiable densities:
$$ \frac{1}{J(X+Y)-J_{XY}} \geq \frac{1}{J_{XX}-J_{XY}} + 
\frac{1}{J_{YY}-J_{XY}}.$$
Equality holds when $X,Y$ are multivariate normal.
\end{proposition}
\begin{proof} Using the conditional representation, Lemma
\ref{lem:fishcondsum}, for any $a,b$:
\begin{eqnarray*}
0 & \leq & \ep \left( a \ra{\rho}{1}(X,Y) + b \ra{\rho}{2}(X,Y) - 
(a+b) \wir(X+Y) \right)^2
\\
& = & a^2 J_{XX} + 2ab J_{XY} + b^2 J_{YY} - (a+b)^2 J(X+Y).\end{eqnarray*}
Now, motivated by the choice of $a,b$ that give equality in the Gaussian case,
we take $a=J_{YY} - J_{XY}$, $b=J_{XX} - J_{XY}$, and rearranging, we obtain 
that:
$$ J(X+Y) \leq \frac{ J_{XX} J_{YY} - J^2_{XY} }{J_{XX} + J_{YY} - 2J_{XY}},$$
and subtracting $J_{XY}$ from both sides we obtain the result.
\end{proof}
\begin{lemma}
If $(X_t,Y_t) = \vc{X}+\vc{Z}_{Ct}$, where
$\vc{Z}_{C t} \sim N(0,C t)$, and $W_t = X_t + Y_t$ then
writing $a = J_{XX}-J_{XY}$, $b=J_{YY}-J_{XY}$:
$$ \frac{\partial}{\partial t} 
\left( 2H(X_t,Y_t) - 2H(W_t) \right) \geq \frac{a^2 C_{11} - 2ab C_{12} 
+ b^2 C_{22}}{a+b} \geq 0.$$ \end{lemma}
\begin{proof} 
Johnson and Suhov \cite{johnson3} prove the multivariate de Bruijn identity:
$$\frac{\partial H}{\partial t}(\vc{X}_t) = \frac{1}{2}
\sum_{i,j} C_{ij} J_{ij}(\vc{X} + \vc{Z}_{Ct}), $$
where $J$ is the Fisher matrix $\ep \rho^T \rho$, with $\rho$, the score
vector equal to $\nabla f/f$.
By Proposition \ref{prop:fishdec}
we deduce that:
\begin{eqnarray*}
\lefteqn{\frac{\partial}{\partial t} 
\left( 2H(X_t,Y_t) - 2H(W_t) \right)} \\ 
& = &  C_{11} J_{XX} + 2C_{12} J_{XY} + C_{22} J_{YY} - (C_{11} + 2C_{12} 
+ C_{22}) J(W_t)  \\
& \geq & C_{11} J_{XX} + 2C_{12} J_{XY} + C_{22} J_{YY} \\
& & - (C_{11} + 2C_{12} + C_{22}) \left( 
\frac{J_{XX}J_{YY} - J_{XY}^2}{J_{XX}+J_{YY} - 2J_{XY}} \right)
\end{eqnarray*}\end{proof}

Now, for functions $f(t),g(t)$, we can define $(X_t,Y_t) = (X,Y) + (Z_1,Z_2)$,
where $Z_1$,$Z_2$ are independent, with $Z_1 \sim N(0,f(t))$, 
$Z_2 \sim N(0,g(t))$. We write $v_{X_t}$, $p_{X_t}$ and $\rho_{X_t}$ for 
the variance, density and score function of $X_t$.
This perturbation ensures that densities are smooth
and allows us to use the $2$-dimensional version of the de Bruijn identity:

\begin{condition} \label{cond:posprod}
For all $t$, $\ep \rho_{X_t}(X_t) \rho_{Y_t}(Y_t) \geq 0$.
\end{condition}
Compare Condition \ref{cond:posprod} with Takano's condition \cite{takano}, 
involving the same term:
\begin{condition} \label{cond:takano}
For all $t$,
$\ep \rho_{X_t}(X_t) \rho_{Y_t}(Y_t) \geq \ep M^2_{\lambda,
\lambda^{-1}}$, where $\lambda = \sqrt{\frac{J(X_t)}{J(Y_t)}}$.
\end{condition}
Takano shows that  Condition \ref{cond:takano} implies that the original 
Entropy Power Inequality, Theorem \ref{thm:epi}, holds. With our weaker 
condition, we provide a weaker, though still interesting, result.
\begin{theorem}[Conditional Entropy Power Inequality] \label{thm:cepi}
If Condition \ref{cond:posprod} \\ holds then:
$$ 2^{2H(X+Y)} \geq 2^{2H(X|Y)} + 2^{2H(Y|X)}.$$ \end{theorem}
\begin{proof} 
Taking $f,g$ defined by
$f'=2^{2H(X_t|Y_t)}$, $g'=2^{2H(Y_t|X_t)}$ and 
defining $s(t) = (2^{2H(X_t|Y_t)}  +
2^{2H(Y_t|X_t)})/2^{2H(W_t)}$, 
\begin{eqnarray*}
\lefteqn{s'(t)} \\ 
& \geq & \frac{1}{2^{2H(W_t)}} \left( \left( 2^{2H(X_t|Y_t)} + 2^{2H(Y_t|X_t)} 
\right) \frac{A^2 f' + B^2 g'}{A+B} - f'g' (J(X_t) + J(Y_t)) \right) \\
& \geq & \frac{1}{2^{2H(W_t)}} \left( \frac{(Af'-Bg')^2}{A+B} 
+ f'g'(A+B -J(X_t) - J(Y_t) \right) \\
& = & \frac{1}{2^{2H(W_t)}} \left( \frac{(Af'-Bg')^2}{A+B} 
+ f'g' \ep M_{1,-1}^2 + 2 f'g' \ep \rho_{X_t} \rho_{Y_t} \right) \geq 0,
\end{eqnarray*}
since
$0 \leq \ep M_{1,-1}^2 = J_{XX} - 2J_{XY} + J_{YY} - J(X_t) -J(Y_t) 
- 2 \ep \rho_{X_t} \rho_{Y_t}$.
Hence $s(t)$ is
an increasing function of $t$. Now as $t \tends$, $s(t) \rightarrow 1$,
since $(X,Y)$ tends to an independent pair of normals. Hence $s(0) \leq
1$ and the result follows. \end{proof}
Cover and Zhang \cite{cover2} provide a bound on the entropy $H(X+Y)$,
under the condition that $X$ and $Y$ have the same marginal density $f$. They
show that $H(X+Y) \leq H(2X)$ if and only if 
$f$ is log-concave (that is, the score function is decreasing). Notice that
our Condition \ref{cond:posprod} holds if $X,Y$ are FKG variables with
log-concave densities.

We write $\psi(X, Y) = 
\sup_{x,y} | p_{X,Y}(x,y)/p_X(x) p_Y(y) - 1 |$, the so-called
$\psi$-mixing coefficient. Note that since
\begin{eqnarray*}
\lefteqn{\ep \rho_{X_t}(X_t) \rho_{Y_t}(Y_t) - \frac{\cov(X_t, Y_t)}{v_{X_t} 
v_{Y_t}}} \\
& = & \int \left( \rho_{X_t}(x) \rho_{Y_t}(y) - \frac{xy}{v_{X_t} v_{Y_t}}
\right) \left( p_{X_t,Y_t}(x,y) - p_{X_t}(x) p_{Y_t}(y) \right) dx dy \\
& \geq & - \psi(X_t, Y_t) \sqrt{J(X_t) J(Y_t) - (v_{X_t} v_{Y_t})^{-1}}, \\
\end{eqnarray*}
Condition \ref{cond:posprod} will hold if for all $t$:
$$ \cov(X,Y) = \cov(X_t, Y_t) \geq v_{X_t} v_{Y_t} \psi(X_t, Y_t)
\sqrt{ v_{X_t} J(X_t) v_{Y_t} J(Y_t) -1}.$$
Now, by Lemma \ref{lem:fishinv}
$J(X_t) \leq 1/(J(X)^{-1} + f(t))$, we know that
if $X,Y$ have the same marginals (and wlog variance 1),
then $v_{X_t} J(X_t) v_{Y_t} J(Y_t) -1 \leq (J^2(X) - 1)/(1+f(t)J(X))$.
Thus, we require that:
$$ \cov(X,Y) \frac{\sqrt{1+f(t)J(X)}}{1+f(t)} \geq \psi(X_t, Y_t) 
(J^2(X) - 1).$$
From $f =0$, we deduce that we need $J(X) \leq \sqrt{\cov(X,Y)/\psi(X,Y) + 1}$,
and if $\lim_{t \tends} f(t)^{1/2} \psi(X_t, Y_t) = 0$, we are through.

Although we know that 
$\psi(X_t, Y_t) \leq \psi(X,Y)$, we need some theory of convexity of
mixing coefficients to provide the most natural conditions.


\begin{thebibliography}{1}

\bibitem{barron}
A.R. Barron.
\newblock Entropy and the {Central Limit Theorem}.
\newblock {\em Annals of Probability}, 14:336--342, 1986.

\bibitem{blachman}
N.M. Blachman.
\newblock The convolution inequality for entropy powers.
\newblock {\em IEEE Transactions on Information Theory}, 11:267--271, 1965.

\bibitem{cover2}
T.M. Cover and Z.~Zhang.
\newblock On the maximum entropy of the sum of two dependent random variables.
\newblock {\em IEEE Transactions on Information Theory}, 40:1244--1246, 1994.

\bibitem{johnson3}
O.T. Johnson and Y.M. Suhov.
\newblock Entropy and random vectors.
\newblock {\em Journal of Statistical Physics}, 104:147--167, 2001.

\bibitem{shannon}
C.E. Shannon and W.W. Weaver.
\newblock {\em A Mathematical Theory of Communication}.
\newblock University of Illinois Press, Urbana, IL, 1949.

\bibitem{stam}
A.J. Stam.
\newblock Some inequalities satisfied by the quantities of information of
  {Fisher} and {Shannon}.
\newblock {\em Information and Control}, 2:101--112, 1959.

\bibitem{takano}
S.~Takano.
\newblock The inequalities of {Fisher Information} and {Entropy Power} for
  dependent variables.
\newblock In S.~Watanabe, M.~Fukushima, Yu.V Prohorov, and A.N. Shiryaev,
  editors, {\em Proceedings of the 7th Japan-Russia Symposium on Probability
  Theory and Mathematical Statistics, Tokyo 26-30 July 1995}, pages 460--470,
  Singapore, 1996. World Scientific.

\end{thebibliography}
\end{document}